\documentclass[journal ]{new-aiaa}
\usepackage[utf8]{inputenc}
\usepackage{textcomp}

\usepackage{graphicx,subcaption}
\usepackage{wrapfig}
\usepackage{tikz}
\usetikzlibrary{arrows,decorations.pathreplacing}
\usepackage{amsmath}
\usepackage[version=4]{mhchem}
\usepackage{siunitx}
\usepackage{longtable,tabularx}
\usepackage{float}
\restylefloat{table}
\usepackage{amsfonts}
\usepackage{mathtools} 
\usepackage{tikz}
\usetikzlibrary{arrows,decorations.pathreplacing}
\tikzstyle{underbrace style}=[decorate,decoration={brace,raise=1mm,amplitude=5pt,mirror},color=black]
\tikzstyle{underbrace text style}=[font=\large, below, pos=.5, yshift=-4mm]

\tikzset{
    rect-color/.style = {
        very thick,
        draw = #1,
        fill = #1!40,
        rounded corners = 2pt,
        opacity = 0.5
    }
}

\newtheorem{manualtheoreminner}{Theorem}
\newenvironment{manualtheorem}[1]{%
  \renewcommand\themanualtheoreminner{#1}%
  \manualtheoreminner
}{\endmanualtheoreminner}

\newtheorem{lemma}{Lemma}
\newtheorem{corollary}{Corollary}
\newtheorem{proposition}{Proposition}
\newtheorem{theorem}{Theorem}
\newtheorem{claim*}{Claim}

\newtheorem{question}{Question}
\newtheorem{example}{Example}
\DeclareMathOperator{\cl}{Cl}
\DeclareMathOperator{\pref}{pref}
\DeclareMathOperator{\suff}{suff}
\newcommand{\A}{\mathbb{A}}

\title{Finite and infinite closed-rich words\footnote{The first author is supported by the Ministry of Science and Higher Education of the Russian Federation (agreement no. 075–15–2022–287).\\ The second author is supported by Foundation for the Advancement of Theoretical Physics and Mathematics ``BASIS''.}}

\author{Olga Parshina$^{1}$ and Svetlana Puzynina$^{1,2}$}
\affil{$^1$Saint Petersburg University, 7/9 Universitetskaya nab., St. Petersburg, 199034 Russia\\ $^2$
Sobolev Institute of Mathematics, 4 Acad. Koptyug av., Novosibirsk, 630090 Russia\\
\textnormal{\texttt{\{parolja,s.puzynina\}@gmail.com}}}

\begin{document}

\maketitle

\begin{abstract} {\begin{singlespace} Abstract: \textnormal { A word is called closed if it has a prefix which is also its suffix and there is no internal occurrences of this prefix in the word. In this paper we study words that are rich in closed factors, i.e., which contain the maximal possible number of distinct closed factors. As the main result, we show that for finite words the asymptotics of the maximal number of distinct closed factors in a word  of length $n$ is $\frac{n^2}{6}$. For infinite words, we show there exist words such that each their factor of length $n$ contains a quadratic number of distinct closed factors, with uniformly bounded constant; we call such words infinite closed-rich. We provide several necessary and some sufficient conditions for a word to be infinite closed rich. For example, we show that all linearly recurrent words are closed-rich. We provide a characterization of rich words among Sturmian words. Certain examples we provide involve non-constructive methods.} \end{singlespace}}

\end{abstract}
\section{Introduction}

\vspace{0.2cm}
Various questions that concern counting factors of a specific form in a word of length $n$ have been studied in combinatorics on words. 
Several studies have been devoted to the words that are extremal with respect to the proportion of factors with a given property. For example, an extensive study has been performed on the problem of counting the maximal repetitions (runs) in a word of length $n$. It has been shown in \cite{DBLP:conf/focs/KolpakovK99} that the maximal number of runs in a word is linear, and it was conjectured to be $n$. Subsequently, there was a lot of research performed to find the bound~\cite{DBLP:journals/tcs/CrochemoreIT11}.
Recently, the conjecture has been proved with a remarkably simple argument, considering numerous attempts to solve it~\cite{DBLP:journals/siamcomp/BannaiIINTT17}. 
We remark that questions about counting regular factors in a word are often non-trivial. For example, the problem of bounding the number of distinct squares in a string: A.S.~Fraenkel and J.~Simpson showed in 1998 \cite{DBLP:journals/jct/FraenkelS98} that a string of length $n$  contains at most $2n$  distinct squares, and conjectured that the bound is actually $n$. After several improvements, the bound of $\frac{11}{6}n$ was proved in \cite{DBLP:journals/dam/DezaFT15}, but the conjecture remains unsolved. 

A related problem concerns counting palindromic factors. It is easy to see that a word of length $n$ can contain at most $n+1$ distinct palindromes  (see e.g.~\cite{DROUBAY2001539}). 
Such words are called \emph{rich in palindromes}, and there also exist infinite  words such that all their factors are rich. Words rich in palindromes have been characterized in~\cite{DBLP:journals/ejc/GlenJWZ09}. 
Words containing few palindromes were studied in~\cite{BrlekPalindromes,FICI20131}. 
Recently some related questions about counting generalizations of palindromes have been studied, e.g.\ privileged factors \cite{DBLP:journals/tcs/Peltomaki13} and $k$-abelian palindromes \cite{DBLP:journals/iandc/CassaigneKP18}.
 
We are interested in counting the factors that are called \emph{closed}. A finite word is called \emph{closed} if it has length $\leq1$ or it is a complete first return to some proper factor, i.e.\ it starts and ends with the same word that has no other occurrences but these two. Otherwise the word is called \emph{open}.
The terminology closed and open was introduced by G.~Fici in~\cite{Fici_oc}; 
for more information on closed words see~\cite{surveyFici}. The notion of a closed word is actually the same as the notion of complete return word. The name return word is usually referred to factors of an infinite word and
is used to study its properties. It can be regarded as  a discrete analogue of the first return map in dynamical systems. For example, F.~Durand characterized primitive substitutive words using the notion of a return word \cite{DBLP:journals/dm/Durand98}. Return words also provide a nice characterization of the family of Sturmian words \cite{DBLP:journals/ejc/Vuillon01}. 
The explicit formulae for the functions of closed and open complexities for the family of Arnoux--Rauzy words, which encompass Sturmian words, were obtained in~\cite{ParZamb}. 
In~\cite{ParPost2020}, the authors proved a refinement of the Morse--Hedlund theorem (see~\cite{MorseHed}) providing a criterion of periodicity of an infinite word in terms of closed and open complexities.

The concept of a closed factor has recently found applications in string algorithms.
The \emph{longest closed factor array} (LCF array) of a string $x$ stores for every suffix of $x$ the length of its longest closed prefix. It was introduced in~\cite{BADKOBEH201623} in connection with closed factorizations of a string. Among other things, the authors presented algorithms for the factorization of a given string into a sequence of longest closed factors and for computing the longest closed factor starting at every position in the string. In~\cite{Bannaietal2015}, the authors presented the algorithm of reconstructing a string from its LCF array. See also~\cite{alamro2020efficient} for some generalizations.

\vspace{0.5mm}
It is easy to show that each word of length $n$ contains at least $n+1$ distinct closed factors~\cite{Fici_low}. In this paper, we study closed-rich words, i.e., words containing the maximal number of distinct closed factors among words of the same length. 
We show the asymptotics of $\sim\frac{n^2}{6}$ for this number. The family of words achieving this bound is constructed in  Proposition~\ref{pr:lowerbound}, and the proof for the upper bound, which is quite technical, is provided in Theorem~\ref{th:main1}. 
We also extend the notion of closed-rich words to infinite words, requiring each factor to contain a quadratic number of distinct closed factors, with uniformly bounded constant; we call such words infinite closed-rich. 
We show that infinite closed-rich words exist, providing a series of families of such words. We prove several necessary and some sufficient conditions for a word to be infinite closed rich. For example, we show that all linearly recurrent words are closed-rich (Proposition~\ref{pr:infrich_linrec}); we provide a characterization of rich words among Sturmian words (Corollary~\ref{cor:Sturmian}). Some examples are non-constructive: we make use of a result by Beck~\cite{Beck}, showing the existence of a word with certain properties using Lov\'asz local lemma.

Some initial results on the topic have been reported at CSR~2021~\cite{DBLP:conf/csr/ParshinaP21}. In particular, an upper bound for the number of distinct closed  factors in a word of length $n$ has been obtained and the asymptotics of $\sim\frac{n^2}{6}$ has  been conjectured. In this paper we prove the conjectured asymptotics. 
In the case of infinite words, we generalize results from \cite{DBLP:conf/csr/ParshinaP21}: in particular, we make an observation that linearly recurrent words are closed-rich, which leads to a nicer and less technical characterization of closed-rich words in the family of Sturmian words. Besides that, the proofs making use of non-constructive words are also new compared to the conference version.


\section{Preliminaries}

In this section, we introduce the necessary definitions and notation, and prove some initial properties of closed-rich words.

Let $\A$ be a finite set called an alphabet. A finite or an infinite word $w=w_0 w_1 \cdots$ on $\A$ is a finite or infinite sequence of symbols from $\A$. 
For a finite word $w=w_0\cdots w_{n-1}$, its \emph{length} is $|w|=n$. 
We let $\varepsilon$ denote the empty word, and we set $|\varepsilon|=0$. 
A word $v$ is a \emph{factor} of a finite or an infinite word $w$ if there exist words $u$ and $y$ such that $w$ can be represented as their concatenation $w=uvy$. 
If $u=\varepsilon$, then $v$ is a prefix, and if $y=\varepsilon$, then $v$ is a suffix of $w$. 
If a finite word $w$ has a proper prefix $v$ which is also its suffix, then $v$ is called a \emph{border} of $w$. 
If the longest border of a word $w$ occurs in $w$ only twice (as a prefix and as a suffix), then $w$ is \emph{closed}. By convention, if $w$ is the empty word or a letter, then it is closed. 

It is not hard to see that a word of length $n$ contains at least $n+1$ distinct closed factors; G.~Fici and Z.~Lipt\'ak characterized words having exactly $n+1$ closed factors~\cite{Fici_low} .
In the same paper they showed that there are words containing $\Theta(n^2)$ many distinct closed factors. The example they provided is a binary word with $\sim \frac{n^2}{32}$ closed factors.
We say that a finite word $w$ is \emph{closed-rich} if it contains at least as many distinct closed factors as any other word of the same length and on the alphabet of the same cardinality. 

If there exists an integer $t$ such that for each $i$ ($i< |w|-t$ in the case $w$ is finite) the equality $w_{i+t}=w_i$ holds, then $t$ is called a \emph{period} of $w$. 
Let $s= \frac{|w|}{t}$ and let $u$ be the prefix of $w$ of length $t$. We say that $w$ has \emph{exponent} $s$ and write $w=u^s$. 
The notation $w=u^{k+}$ means that $w$ has exponent $s>k$ for an integer $k$.
The word $u$ is called the \emph{fractional root} of $w$.
The word $w$ is \emph{primitive} if its only integer exponent is 1.
Hereinafter we always assume $t$ to be the shortest period of $w$, and thus, $s$ to be the largest exponent of $w$. 

The following properties follow directly from the definitions. 

\begin{proposition}\label{pr:exp2closed}
Any word with exponent at least two is closed.
\end{proposition}

\begin{proposition}\label{le:cube} Let $w$ be a word of exponent $3$ and of length $n$. Then all its factors of length at least $\frac{2n}{3}$ are closed, and moreover, all of them except for the factor of length $\frac{2n}{3}$ are unioccurrent.
\end{proposition}

De~Bruijn graph of order $n$ on an alphabet $\A$ is the directed graph whose set of vertices (resp. edges) consists of all words over $\A$ of length $n$ (resp. $n+1).$ There is a directed edge from $u$ to $v$ labeled $w$ if $u$ is a prefix of $w$ and $v$ a suffix of $w.$
We call a Hamiltonian path in this graph a \emph{de~Bruijn word}.

\begin{proposition} \label{pr:cube_dB} 
Let $n=3\cdot |\A|^k$ for an integer $k$, and $v$ be a de~Bruijn word of length $\frac{n}{3}$. Then $w=v^3$ has $\sim \frac{n^2}{6}$ distinct closed factors. 
\end{proposition}

\textit{Proof.}
Due to Proposition~\ref{le:cube}, all factors that are longer than $\frac{2n}{3}$ are closed and distinct (there are $\sim\frac{n^2}{18}$ of those). 
All words of length $\frac{n}{3}+\log \frac{n}{3} \leq l \leq \frac{2n}{3}$ are also closed with corresponding border of length $l-\frac{n}{3}$ (there are $\sim \left(\frac{n}{3}\right)^2$ of distinct factors of these lengths).

If a factor of length less than $\frac{n}{3}+\log\frac{n}{3}$ is closed, then its border is shorter than $\log\frac{n}{3}$, because all factors of de~Bruijn word of length at least $\log\frac{n}{3}$ are unioccurrent. Thus, there are not more than $\frac{n}{3} \cdot\log\frac{n}{3}$ closed factors that are shorter than $\frac{n}{3}+\log\frac{n}{3}$. 
$\square$

The construction from the previous proposition gives only words of length $n=3\cdot |\A|^{k}, k\geq0$, but it could be easily modified to other lengths (see Proposition~\ref{pr:lowerbound}).

The following example shows that if we change one letter in the middle of a word, the total number of closed factors can change dramatically:

\begin{example} {\normalfont 
Let us show that the number of closed factors can change from linear to quadratic when changing only one letter in a word. It is easy to see that the word $a^nba^nba^n$ has quadratic number of closed factors. After replacing the leftmost occurrence of $b$ with $a$, we obtain the word $a^n a a^nba^n$ with linear number of closed factors.}
\end{example}
For a finite word $w$, we let $\cl(w)$ denote the number of distinct closed factors of $w$. 
We now provide a trivial upper bound on $\cl(w)$.


\begin{proposition}\label{pr:upper}
For each word $w$ of length $n\geq 7$, one has $\cl(w)\leq \frac {n^2}{4}$.
\end{proposition}

\textit{Proof.} For a word $u$ and a letter $a$, let us denote by $t$ the longest repeated suffix of $ua$, and by $z$ the longest repeated suffix of $t$.  
Clearly, the number of new closed factors ending in the last letter of $ua$ is at most  $|t|-|z|\leq  \left\lfloor \frac{|ua|}{2} \right\rfloor $ when $u$ is non-empty. For $n=0$ we have one closed factor (the empty word), for $n=1$ we add one closed factor (a letter). So, building $w$ letter by letter, we get at most
$2+\sum_{i=2}^n \lfloor \frac{i}{2} \rfloor = 2+\frac{n(n+1)}{4}-\lceil \frac{n}{2}\rceil$ closed factors in $w$. The claim follows. 
$\square$

In Section \ref{sec:fin} we prove an upper bound of $\sim \frac{n^2}{6}$ and provide the asymptotics of the maximal number of distinct closed factors a word of length $n$ can contain.


Sometimes it will be convenient for us to consider cyclic words. For a word $w$, we can consider a corresponding cyclic word as the class of all its cyclic shifts. By a closed factor of a cyclic word we mean a closed factor of some its shift. 

\begin{lemma}\label{le:cyclic}
Let $u$ be a primitive finite word of length $k$, then its cyclic square has at most $k^2$ distinct closed factors.
\end{lemma}

\textit{Proof.}
Let $\hat u$ be the cyclic square of $u$ with the first letter $\hat u_0= u_0$.
The following observations constitute the proof of the lemma. Basically, we count (left) borders giving rise to distinct closed words.

Each occurrence of a factor of $\hat u$ is a (left) border of at most one closed factor. Since $\hat u$ is a cyclic square, in order to count  borders giving rise to distinct closed words, we can only consider the factors of $\hat u$ starting in its first half $\hat u_0\cdots \hat u_{k-1}=u_0\cdots u_{k-1}$ (borders starting in the second half of  $\hat u$ give rise to the same closed words). The border of a closed factor of $\hat u$ cannot be longer than $k$, otherwise $u$ is not a primitive word.
The number of factors of $\hat u$ that start in its first half and are not longer than $k$ is $k^2$.
The statement follows.
$\square$

\begin{lemma}\label{le:exp3} Let $w$ be a word with exponent at least $3$. 
Then the number of closed factors in the word is at most $k^2+{(n-3k)}{k}+\frac{1}{2}(k+1)k$. \end{lemma}

\textit{Proof.} We count closed factors by their lengths: $k^2$ is the upper bound on the number of closed factors of length at most $2k$ (given by the bound for a cyclic word of length $2k$ from Lemma~\ref{le:cyclic}); $(n-3k)k$ counts $k$ closed factors of each length $2k+1, \dots, n-k$; and $\frac{1}{2}(k+1)k$ counts long ones as the sum of an arithmetic progression ($k$ factors of length $n-k+1$, $k-1$ of length $n-k+1$, \dots, and $1$ of length $n$).  $\square$

\begin{corollary} Words of length $n$ and of exponent greater than $3$ have not more than $\frac{n^2}{6}$ factors asymptotically.\end{corollary}

 \textit{Proof.} We estimate the bound from Lemma~\ref{le:cube}: for words with exponent $t \geq 3$ we count closed factors and get the function $c(t)=n^2(\frac{1}{t}- \frac{3}{2t^2})$. The maximum value of $c(t)$ is $\frac{n^2}{6}$ and it is achieved for $t=3$. $\square$

We say that a finite word $w'$ is a \emph{cyclic shift} (or a conjugate) of a finite word $w$ if there exist words $u$ and $v$ such that $w=uv$ and $w'=vu$.

\begin{lemma}\label{le:cyclicube}
Let $w$ be a word of exponent $\alpha \geq 3$ and $v$ its primitive root, so that $w=v^{\alpha}$. Then for any cyclic shift $v'$ of $v$, the word $v'^{\alpha}$ contains the same number of closed factors as $w$.
\end{lemma}

\textit{Proof.} The statement follows from the counting in the proof of Lemma~\ref{le:exp3}: the set of closed factors consists of all long factors (their numbers are clearly the same since we only take lengths into account) and short closed factors which are also closed factors of the cyclic square, so their sets are the same. 
$\square$

The following example shows that the statement of Lemma~\ref{le:cyclicube} does not hold for squares. Moreover, it is possible that taking a cyclic shift of the root changes the number of closed factors in a word from linear to quadratic.

\begin{example} {\normalfont  The number of closed factors in $a^{\frac{n}{2}}ba^nba^{\frac{n}{2}}$ is quadratic, while in its cyclic shift $a^n b a^n b$ it is linear.}
\end{example}

\section{Infinite rich words}

In the previous section we provided some examples of finite words containing a quadratic number of distinct closed factors. In this section, we extend this property to infinite words. Namely, we consider infinite words for which there exists a constant $C$ such that for each $n\in\mathbb{N}$ each factor of $w$ of length $n$ contains at least $Cn^2$ distinct closed factors. We call such words \emph{infinite closed-rich}. We remark that in this definition we do not require the constant to be optimal; however, a natural question is optimizing this constant (see Question~\ref{qu:absolute_rich}). We show that infinite closed-rich words exist, and provide some families of examples of words that are infinite closed-rich, as well as families of words that are not.

An \emph{exponent} (or index) of an infinite word is defined as the supremum of exponents of its factors. 

\begin{proposition}\label{pr:infrich_powers} Let $w$ be an infinite word of unbounded exponent. Then $w$ is not closed-rich.
\end{proposition}
\textit{Proof.} Let us show that a word $v$ of length $n$ and of exponent $k$ contains at most $\frac{n^2}{k}$ closed factors. 
To estimate the number of distinct closed factors in $v$ from above, it is enough to estimate the number of all closed factors that begin inside the first period. Indeed, all closed factors beginning later are equal to factors beginning in the first period. The length of the period of $v$ is $l=\frac{n}{k}$. 
In each position $1, \dots, l$ there are less than $n$ possible borders, hence the total number of distinct closed factors in $v$ is at most $nl=\frac{n^2}{k}$.
Since $w$ is of unbounded exponent, $k$ can be any number, so there is no uniform constant $C$ from the definition of infinite closed-rich words, and hence the word $w$ is not closed-rich. 

\bigskip

An infinite word $w$ is called \emph{linearly recurrent} if there exists a constant $K$ such that the length of each return word to each factor $u$ of $w$ is at most $K|u|$. 

\begin{proposition}\label{pr:infrich_linrec} Let $w$ be an infinite linearly recurrent aperiodic word. Then $w$ is  closed-rich.
\end{proposition}

\textit{Proof.} Since $w$ is linearly recurrent, there exists a constant $K$ such that the length of each return word to a factor $u$ of of $w$ is at most $K|u|$. It means that each factor $v$ of $w$ of length at least $2Kn+n-1$ 
contains all factors of length at most $n$  at least twice. There are at least $\frac{n^2}{2}$ such factors in total (since the word is aperiodic, by Morse and Hedlund theorem~\cite{MorseHed} there are at least $n+1$ factors of each length; and there are $n$ possible lengths); each factor is a border of at least one closed factor, and all these factors are distinct. Thus, we get a quadratic lower bound on the number of distinct closed factors in each factor of $w$, with the constant depending only on the parameter $K$ of linear recurrence. $\square$








As another corollary from Propositions \ref{pr:infrich_powers} and \ref{pr:infrich_linrec}, we get a characterization of closed rich words in the family of Sturmian words. \emph{Sturmian words} are usually defined as infinite words with the smallest possible number of distinct factors of each length among aperiodic words ($n+1$ factor of each length $n\geq 1$). Sturmian words are known to be rich in palindromes~\cite{DBLP:journals/ejc/GlenJWZ09}. 
Sturmian words admit various characterizations, the one we use here is via balance. 
An infinite word is called \emph{balanced} if for each its two factors $u$ and $v$ of the same length and each letter $a$ one has $||u|_a-|v|_a|\leq 1$, where $|u|_a$ denotes the number of occurrences of the letter $a$ in $u$. 
Sturmian words can be equivalently defined as aperiodic binary balanced words~\cite{DBLP:journals/mst/CovenH73}. It is known that balanced binary words have uniform frequencies of letters; for a Sturmian word $s$ the frequency $\alpha$ of 1's in it is irrational and is called a \emph{slope} of $s$. Consider a continued fractional expansion of $\alpha$: 

$$\alpha=m_0 +
\frac{1}{m_1+\frac{1}{m_2+\frac{1}{m_3+\frac{1}{\dots}}}},$$ where
$m_0, m_1, \dots$ are integers.

It is known that a Sturmian word is linearly recurrent if and only if the sequence $(m_i)_{i\geq 1}$ is uniformly bounded, which is also equivalent to the fact that the Sturmian word is of bounded index. For these facts and for more information on Sturmian words we refer the reader to Chapter~2 of~\cite{Lothaire}.

\begin{corollary}\label{cor:Sturmian}
A Sturmian word is closed-rich if and only if the directive sequence of its slope is bounded.
\end{corollary}

A \emph{morphism} $\varphi$ is a map on the set of all finite words on the alphabet $\A$ such that $\varphi(uv) = \varphi(u)\varphi(v)$ for all finite words $u,v$ on $\A$. 
The domain of the morphism $\varphi$ can be naturally extended to infinite words by $\varphi(w_0w_1w_2 \cdots) = \varphi(w_0)\varphi(w_1)\varphi(w_2) \cdots$. A \emph{fixed point} of a morphism $\varphi$ is an infinite word $w$ such that
$\varphi(w) = w$.
A morphism $\varphi$ is \emph{primitive} if there exists a positive integer $l$ such that the letter $a$ occurs in the word $\varphi^l(b)$ for each pair of letters $a, b \in \A$. 
Since words generated by primitive morphisms are linearly recurrent~\cite{DBLP:journals/dm/Durand98}, the following holds.

\begin{corollary} Infinite words generated by primitive morphisms are closed-rich.
\end{corollary}

In particular, it follows that the Thue-Morse and the Fibonacci words are closed-rich.

\bigskip

Summing up the above, we showed that closed-rich words must be of bounded index, and that all linearly recurrent aperiodic words are closed-rich. A natural question is whether there exist uniformly recurrent words of bounded index, but are not closed-rich. Such words exist, and one of the ways to obtain them is using a beautiful probabilistic and non-constructive result due to Beck:

\begin{theorem}[\cite{Beck}]
For any real $\varepsilon
>0$, there exist an integer $N_{\varepsilon}$ and an infinite binary word $w$ such
that for every factor $x$ of $w$ of length $n > N_{\varepsilon}$,
all occurrences of $x$ in $w$ are separated by a distance at least
$(2-\varepsilon)^n$.\end{theorem}

Let us call the words satisfying the theorem \emph{Beck words} with parameter $\varepsilon$, or \emph{$\varepsilon$-Beck words}.

\begin{theorem} \label{th:Beck1} For each $\varepsilon<1$, any Beck word with parameter $\varepsilon$ is not closed-rich.  Moreover, for each positive integer $C$, for sufficiently large length $N$ all its factors of length $N$ contain less than $CN^2$ closed factors.  \end{theorem}

\textit{Proof.}
First, let us notice that for each factor $v$ of length $N>(2-\varepsilon)^n$, $n>N_{\varepsilon}$, of a Beck word, the borders of factors of $v$ have length less than $n$. 
Indeed, since the factors of length $n$ of the Beck word are separated by at least $(2-\varepsilon)^n$ letters, all factors of length at least $n$ in $v$ are unioccurrent, so they cannot be borders of factors of $v$. 
Thus, each position of $v$ can be a beginning of less than $n$ closed words. Therefore, the total number of closed factors of $v$ (even with repetitions) does not exceed $nN<N \log_{2-\varepsilon}N$. The claim follows.
$\square$

For an infinite word $w$, one can define its \emph{shift orbit closure} as the set of infinite words whose sets of factors are included in the set of factors of $w$. 

\begin{corollary} There exists an infinite uniformly recurrent word of bounded exponent which is not closed-rich.  
\end{corollary}
\textit{Proof.} Since in a shift orbit closure of each word there exists a uniformly recurrent word, we can consider a uniformly recurrent word in a shift orbit closure of a Beck word. We claim that it gives the example. 
Indeed, due to Theorem~\ref{th:Beck1}, each long enough factor of this word contains less than quadratic number of closed factors, so it remains to show that any Beck word is of bounded exponent.

Suppose an $\varepsilon$-Beck word $w$ is of unbounded exponent. Two cases are possible: either there exists a factor $x$ such that for each $k$ the word $x^k$ is a factor of $w$, or there exists a series of words $(x_i)_{i\geq 1}$ of increasing lengths such that $x_i^i$ is a factor of $w$. 
In the first case, we get a contradiction considering a power $x^2j$ with  $|x^j|>N_{\varepsilon}$: there are two copies of $x^j$ at distance 0. In the second case, consider $x_i$ with $i>1$ and $|x_i|>N_{\varepsilon}$, then the factor $x_i^2$ gives the contradiction with the definition of Beck word.$\square$








\begin{proposition} \label{pr:infrich} Let $w$ be an infinite word, and let $C>2$, $\alpha<1$ be two constants. If for each $n$ each factor of $w$ of length $n$ contains a factor of exponent at least $C$ and of length of period at least $\alpha n$, then $w$ is infinite closed-rich.\end{proposition}

\textit{Proof.} Let $v$ be a factor of $w$ of length $n$. There is a factor $u$ of $v$ of period $k\geq \alpha n$ and of exponent $C'\geq C>2$, hence its length $l = C'k \geq C \alpha n $. 
To count closed factors of $u$ we use the following two observations. All factors of $u$ of length greater than $l-k$ are distinct.
Each factor of $u$ of length at least $2k$ has exponent at least~2 and hence is closed by Proposition~\ref{pr:exp2closed}.

If $C'\leq 3$, there are at least
$\sum_{j=2k}^l (l-j)\geq\frac{(C'-2)^2}{2}k^2 \geq \frac{(C-2)^2(\alpha n)^2}{2}$ closed factors. 

If $C'>3$, then, in addition to the closed factors longer than $l-k$, the word $u$ has $k$ distinct closed factors of each length between $2k$ and $l-k$.
Thus, there are at least
$(l-3k)k + \sum_{j=l-k}^l (l-j)  \geq (C'-3)k^2 +\frac{k^2}{2} = \frac{C'-5}{2}k^2 > \frac{(\alpha n)^2}{2} = \frac{\alpha^2 n^2}{2}$ closed factors.

Therefore, the constant in the definition of infinite rich words is given by $\min\left(\frac{\alpha^2}{2},  \frac{(C-2)^2\alpha^2}{2}\right)$.  
$\square$

\begin{example} {\normalfont
Let $w$ be a fixed point of the morphism $\varphi: a\to abbba, b\to abbbb$. We show that it is infinite closed-rich. 
Indeed, each block $\varphi^k(c)$ for $c\in\{a,b\}$ has length $5^k$ and contains a cube with the period $5^{k-1}$. Clearly, each factor of length at least $2\cdot 5^k-1$ contains a block $\varphi^k(c)$. The maximal length where we cannot guarantee the occurrence of the block $\varphi^{k+1}(c)$, is $2\cdot 5^{k+1}-2$.
Thus, we can apply Proposition~\ref{pr:infrich} with $C=3$ and $\alpha = \frac{1}{(2\cdot 5^2)}=0.02$. 
We remark that this word can also be seen as a Toeplitz word~\cite{toeplitzdef} with pattern $baaa?$, and that this construction can be easily generalized to other morphic and Toeplitz words.}
\end{example}

\section{Finite closed-rich words}\label{sec:fin}
The main goal of this section is to prove Theorem~\ref{th:main} providing an upper bound on the number of closed factors that a finite word can contain. Together with the lower bound from Proposition~\ref{pr:lowerbound}, it gives an asymptotics of $\sim \frac{n^2}{6}$.

The set of distinct closed factors of a word $w$ can naturally be split into two sets, the set of closed words of length at least 2, which have non-empty borders, and the set of closed words of length at most 1, i.e., letters and the empty word. We let $\cl'(w)$ denote the number of closed words of length at least $2$ (``long'' closed factors), and $\cl^0(w)$ denote the number of ``short'' closed factors, so that $\cl(w)=\cl'(w)+\cl^0(w)$.

\begin{theorem}\label{th:main}
For a finite word $w$ of length $n$, the following holds: 
\[\cl'(w)\leq  \frac{n^2}{6}+\frac{n}{6}.\]
\end{theorem}
Clearly, since $\cl(w)=\cl'(w)+|\A|+1\leq \cl'(w)+ n+1$, we can rewrite the statement of Theorem~\ref{th:main} as follows.

\begin{manualtheorem}{\ref{th:main}$'$}\label{th:main1}
For a finite word $w$ of length $n$, the following holds:
\[\cl(w)\leq \frac{1}{6}n^2+\frac{7}{6}n+1.\]
\end{manualtheorem}
For a finite word $w$ of length $n$, we denote its prefix of length $n-1$ by $w^-$ and its suffix of length $n-1$ by $^-w$.
The following lemma constitutes the key part of the proof of Theorem~\ref{th:main}.

\begin{lemma}\label{upperC}
Let $w$ be a word of length $n$. If $\cl'(w)-\cl'(w^-) > \frac{n}{3}$ and $\cl'(w)-\cl'(^-w) > \frac{n}{3}$, then $\cl'(w) \leq \frac{1}{6} n^2 + \frac{1}{6}n$.
\end{lemma}

\textit{Proof.}
In the proofs of the lemma and of Theorem~\ref{th:main}, we only talk about long closed factors (of length at least 2).
Let $t$ denote the longest repeated suffix of $w$, $z$ denote the longest repeated suffix of $t$, and $c$ be the letter preceding the last occurrence of $z$ in $w$, so that $cz$ is the shortest unrepeated suffix of $t$. Clearly,  \begin{equation}\label{eq:cl}\cl'(w)-\cl'(w^-)\leq |t|-|z|.\end{equation} 
Two cases are possible: the last and the penultimate occurrences of $t$ in $w$ might intersect, or not.
If they intersect, we have a power as a suffix of $w$; let $l$ be the period of this power. In this case we denote the suffix of $t$ of length $l$ by $x$. In this case $z$ is of length at least $|t|-l$; thus $|t|-|z|$ is at most $l$, and thus $l>\frac{n}{3}$. So these intersecting borders $t$ form a word of exponent less than 3 with primitive root $x$.

If they do not intersect, we set $x=t$. We choose $C> \frac{1}{3}$ to be a constant such that $l-|z|=Cn$. The same reasoning gives the prefix $t'$ with its longest repeated prefix $z'$, the inequality $\cl'(w)-\cl'(^-w)\leq |t'|-|z'|$, the constant $C'=\frac{l-|z'|}{n}>\frac{1}{3}$ and the word $x'$ of length $l'$.

There are several possibilities of how the word $w$ can look like depending on whether the occurrences of $t$ intersect or not and whether the penultimate occurrence of $t$ starts from the beginning of the word or not.

Cases:

\begin{enumerate}
    \item $w=xvx$ for a non-empty word $v$; 
    \item $w=ux^{s}= x'^{s'} u'$ for words $u,u'$ and  $2\leq s,s'<3$;
    \item $w=ux^s=x'v'x'u'$ for non-empty words $u,u',v$;
    \item $w=uxvx=x'v'x'u'$ for some non-empty words $v,v'$ and words $u,u'$.
\end{enumerate}

\bigskip

\begin{figure}[H] 
\centering
\begin{tikzpicture}[thick,scale=0.9, every node/.style={scale=0.75}]
\centering
\pgfmathsetmacro\l{16}
\pgfmathsetmacro\i{1}
\pgfmathsetmacro\u{6}
\pgfmathsetmacro\z{1.5}
\Large
\draw (0,0) -- (\l,0);

\draw [underbrace style] (0,0) -- (\u,0) node [underbrace text style] { $l$};
\draw [underbrace style]  (\l-\u,0) -- (\l-\z,0) node [underbrace text style] { $Cn > \frac{n}{3}$};

\draw (0,0) arc (180:0:{\u/ 2}  and 0.8)
node[midway, above] {$x$};
\draw (\l-\u,0) arc (180:0:{\u/ 2}  and 0.8)
node[midway, above] {$x$};
\draw (\u,0) arc (180:0:{(\l-2*\u)/ 2}  and 0.75)
node[midway, above] {$v$};

measures
\foreach \x in {0,\u-\z,\u-\z-\z/5, \u,\l-\u,\l-\z,\l-\z-\z/5,\l} 
   \draw(\x,-0.1) -- (\x ,0.1);

\node[above] at (\u-2*\z/3,-0.1) {$z$};   
\node[above] at (\l-2*\z/3,-0.1) {$z$};   
\node[above] at (\u-\z-\z/10,-0.1) {$c$};   
\node[above] at (\l-\z-\z/10,-0.1) {$c$};   
 
\end{tikzpicture}
\caption{The case $w=xvx, |w|=n$.}
    \label{fig:xvx}
\end{figure}
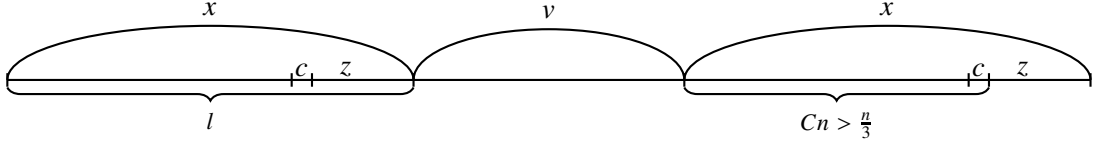

\noindent
1)\ \ The case $w=xvx$ is shown on Fig.~\ref{fig:xvx}.

Let us count borders that are suffixes of the rightmost occurrences of closed factors of $w$. So, for each position $j$, we consider the longest factor $x''$ ending at $j$ which occurs in $w$ earlier, and consider its longest repeated suffix $z''$. If we denote by $C_j$ the number of rightmost occurrences of closed factors of $w$ ending in $j$, we have $C_j = |x''|-|z''|$.

Let us note that if a factor ending at the position $j\in[1,Cn]$ is a border of some closed factor, than it is not its rightmost occurrence, and we do not need to count it, so $C_j=0$.

Every border that ends at the position $j\in[Cn+1,n]$ cannot contain the first occurrence of $cz$, because it would add an extra occurrence of $cz$ in the prefix $x$ of $w$. 
Thus, every border that ends at the position $j\in[Cn+1,2Cn]$ cannot be longer than $j-Cn$, thus cannot be a border of more than $j-Cn$ closed factors, i.e.\ $C_j\leq j-Cn$. 


Now consider $j\in[2Cn+1,n]$. We now show that in this case $C_j\leq Cn$.
Let the starting position of the longest right border that ends in $j$ be $t_2$; we have $t_2\in[Cn+1, n]$. Let us denote by $t_1$ the position of the occurrence of this border as prefix of the corresponding closed factor. Let us consider its left border. It can either contain the first occurrence of $cz$ in $w$, or it ends before the position $l=|x|$, or it starts after the position $Cn$.

In the first case we need to be sure that $cz$ does not occur in the suffix $x$ of $w$ as a proper factor. Thus we need $t_2 + Cn -t_1 < n - l$, i.e. $t_2-t_1 < n-Cn-l < \frac{n}{3} < Cn$. 
In the second case the left border $x''$ of this closed word is a factor of $x$. We make use of the following claim:

\medskip

\textbf{Claim. }\textit{Each factor $y$ of $x$ of length $|y|> Cn$ has a repeated suffix of length at least $|y|-Cn$.} 

\emph{Proof}. 
Let the first and the last positions of $y$ in $x$ be $y_1$ and $y_2$ respectively, and $\hat{z}$ be the longest repeated suffix of $y$.
Let $z_1$ and $z_2$ be the initial and the last positions of penultimate occurrence of $z$ in $x$ respectively. If $z_1 \geq y_1$, then the word $\hat z=x_{Cn+1}\cdots x_{y_1+|y|}$ is the longest repeated suffix of $y$, and $|y|-|\hat z| \leq y_2-|\hat z| = Cn$.

If $z_1<y_1$, then $z=x_{z_1}\cdots x_{z_2}=x_{Cn+1}\cdots x_l$. Cutting the prefix of length $z_1-y_1$ and the suffix of length $l-y_2$ from $z$, we obtain the border and hence a repeated suffix $x_{y_1}\cdots x_{z_2+y_2-l}= x_{Cn+y_1-z_1}\cdots x_{y_2}$ of $y$; its length is  $y_2-Cn-y_1+z_1=|y|-Cn+z_1$.
Thus, $|\hat z| \geq |y|-Cn+z_1$, and $|y|-|\hat z| \leq Cn-z_1\leq Cn$. The claim is proved.

\medskip

Taking $x''$ as $y$ in the claim, we have in this case $C_j\leq Cn$.

In the third case the left border starts after the position $Cn$, i.e. $t_1>Cn$. If it is longer or equal to $\frac{n}{3}$, then $t_2\leq \frac{2n}{3}$, so $t_2-t_1<\frac{n}{3}$, and thus $C_j< \frac{n}{3}$. If it is shorter than $\frac{n}{3}$, then $C_j< \frac{n}{3}$.
Thus, any such factor starting at $t_2$ cannot be a border of more than $Cn$ closed factors.

We obtain the following formula.
\[
\cl'(w) = \sum\limits_{j=1}^{n} C_j \leq \sum\limits_{j=Cn+1}^{2Cn} (j-Cn) + \sum\limits_{j=2Cn+1}^{n} Cn = \frac{Cn(Cn+1)}{2} + Cn(n-2Cn) = -\frac{3}{2}n^2C^2 +\left(n^2 + \frac{n}{2}\right)C.
\]

When $n$ is fixed, this expression reaches its maximum when when $C=\frac{1}{3}+\frac{1}{6n}$. Since $Cn$ is integer, we can take $C=\frac{1}{3}$ for an upper bound. Thus $\cl'(w) \leq \frac{1}{6}n^2 + \frac{1}{6}n$.

\bigskip

2)\ \ Let $w=ux^{s}=x'^{s}u'$ for $2\leq s <3$ and words $u,u'$. 

First we will prove that in this case we must have $u=u'=\varepsilon$ and $s=s'$. Suppose that at least one of them is non-empty, say, $u'$. 

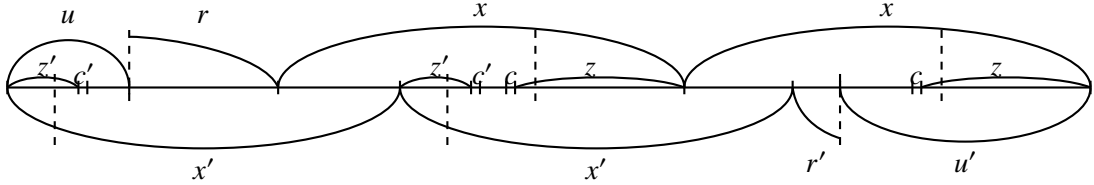
\begin{figure}[H]
\centering
\begin{tikzpicture}[thick,scale=0.9, every node/.style={scale=0.75}]
\centering
\pgfmathsetmacro\l{16}
\pgfmathsetmacro\r{2.2}
\pgfmathsetmacro\u{6} 
\pgfmathsetmacro\v{1.8} 
\pgfmathsetmacro\z{2.8}
\pgfmathsetmacro\uu{5.8} 
\pgfmathsetmacro\vv{3.7} 
\pgfmathsetmacro\z{2.5}
\pgfmathsetmacro\rr{\l-2*\uu-\vv}
\pgfmathsetmacro\zz{\rr+\rr/2}
\Large
\draw (0,0) -- (\l,0);

\draw (0,0) arc (180:0:{\v/ 2}  and 0.7);
\node[above] at (\v/2,0.8) {$u$};

\draw (\v+\r,0) arc (180:0:{\u/ 2}  and 0.9)
node[midway, above] {$x$};
\draw (\l-\u,0) arc (180:0:{\u/ 2}  and 0.9)
node[midway, above] {$x$};
\draw (\v,0.75) .. controls (\v+\r/4,0.75) and (\v+\r-\r/11,0.5)  .. (\v+\r,0);
\node[above] at (\v+\r/2,0.8) {$r$};

\draw (\l-\vv,0) arc (-180:0:{\vv/ 2}  and 0.8)
node[midway, below] {$u'$};
\draw (\uu,0) arc (-180:0:{\uu/ 2}  and 0.9)
node[midway, below] {$x'$};
\draw (0,0) arc (-180:0:{\uu/ 2}  and 0.9)
node[midway, below] {$x'$};
\draw (\l-\vv,-0.75) .. controls (\l-\vv-\rr/4,-0.7) and (\l-\vv-\rr+\rr/11,-0.5)  .. (2*\uu,0);
\node[below] at (\l-\vv-\rr/2,-0.8) {$r'$};

\foreach \x in {0,\v,\v+\r,\v+\r+\u-\z,\v+\r+\u-\z-\zz/8,\l-\u,\l-\z,\l-\z-\zz/8,\l} 
   \draw(\x,-0.1) -- (\x ,0.1);

\draw (\l-\u,0) .. controls (\l-\u-\z/5,0.2) and  (\l-\u-\z+\z/5,0.2)  .. (\l-\u-\z,0);
\node[above] at (\v+\r+\u-\z/1.8,0) {$z$};   
\node[above] at (\v+\r+\u-\z-\zz/18,-0.1) {$c$};  

\draw (\l,0) .. controls (\l-\z/5,0.2) and  (\l-\z+\z/5,0.2)  .. (\l-\z,0);
\node[above] at (\l-\z/1.8,0) {$z$};   
\node[above] at (\l-\z-\zz/15,-0.1) {$c$};   

\foreach \x in {\v,\v+\u,\v+2*\u}
	 \draw[dashed] (\x,0.85) -- (\x ,-0.2);

\foreach \x in {\uu,2*\uu,2*\uu+\rr,\l-\vv, \zz,\zz+\zz/8, \uu+\zz,\uu+\zz+\zz/8} 
   \draw(\x,-0.1) -- (\x ,0.1);

\draw (0,0) .. controls (\zz/5,0.2) and  (\zz-\zz/5,0.2)  .. (\zz,0);
\node[above] at (\zz/1.8,0) {$z'$};   
\node[above] at (\zz+\zz/13,-0.1) {$c'$};  

\draw (\uu,0) .. controls (\uu+\zz/5,0.2) and  (\uu+\zz-\zz/5,0.2)  .. (\uu+\zz,0);
\node[above] at (\uu+\zz/1.8,0) {$z'$};   
\node[above] at (\uu+\zz+\z/15,-0.1) {$c'$};   

\foreach \x in {\rr,\uu+\rr,\rr+2*\uu}
	 \draw[dashed] (\x,-0.85) -- (\x ,0.2);

\end{tikzpicture}
\caption{The case $w=ux^{s}, 2\leq s <3, |w|=n$.}
    \label{fig:4+b}
\end{figure}

Let us consider the shortest unrepeated prefix $z'c'$ of $t'$ (recall that $t'$ is the longest repeated prefix of $w$). We let $r$, $r'$ denote the lengths of fractional parts of exponents of $x$ and $x'$: $r=|x^{s-2}|$,   $r'=|x'^{s'-2}|$ (see Fig.~\ref{fig:4+b}). Clearly, we have $|z'c'|\geq r'$, since otherwise the prefix of length $r'> |z'c'|$ would be repeated in $t'$. Indeed, it occurs in $w$ also at position $l'+1$, which gives its occurrence 
in $x$ at the distance $l'-|u|-r$ from its beginning. 
Let us consider the occurrence of $z'c'$ in the rightmost occurrence of $x$ in $w$.
It cannot occur at the position $2l'+1$, because the letter at the position $2l'+r'+1$ is different from the letter at $l'+r'+1$ by the conditions of the case. It cannot occur either to the left or to the right of this position, because it would add another occurrence of $z'c'$ in $t'$, which is not possible by the definition of $z'c'$. So, $u'$ is empty. A symmetric argument shows that $u$ also must be empty.

It remains to prove that $s=s'$, or, equivalently, that $l=l'$. If $l\neq l'$, then the word $w$ has two periods not greater than $\frac{|w|}{2}$, i.e. $|w|\geq l+l'$. Now we use Fine and Wilf theorem, which says that if a word $w$ has two periods $l$ and $l'$, and its length is at least $l+l'-\gcd(l,l')$, then it also has a period $\gcd(l,l')$ \cite{FW65}. 
In our case it contradicts the assumption that $t=x^{s-1}$ is the longest repeated suffix.

Thus, we have $u,u'=\varepsilon$, $l=l'$ and $l+r-|z|,\ l'+r'-|z'| \geq Cn$. 


\begin{figure}[H] 
\centering
\begin{tikzpicture}[thick,scale=0.9, every node/.style={scale=0.75}]
\centering
\pgfmathsetmacro\l{16}
\pgfmathsetmacro\r{2}
\pgfmathsetmacro\u{7}
\pgfmathsetmacro\z{2.7}
\pgfmathsetmacro\zz{2.5}
\Large
\draw (0,0) -- (\l,0);

\draw (\r,0) arc (180:0:{\u/ 2}  and 1)
 node[midway, above] {$x$};
\draw (\l-\u,0) arc (180:0:{\u/ 2}  and 1)
node[midway, above] {$x$};

\draw (0,0.75) .. controls (\r/4,0.75) and (\r-\r/11,0.5)  .. (\r,0);
\node[above] at (\r/2.7,0.9) {$r$}; 

\foreach \x in {0,\r,\r+\u-\z,\r+\u-\z-\z/9, \r+\u,\l-\u,\l-\z,\l-\z-\z/9,\l,\zz, \zz+\zz/9, \u+\zz, \u+\zz+\zz/9} 
   \draw(\x,-0.1) -- (\x ,0.1);

\draw (\l-\u-\z,0) arc (180:0:{\z/ 2}  and 0.2);

\node[above] at (\r+\u-\z/1.8,0.1) {$z$};   
\node[above] at (\r+\u-\z-\z/20,-0.1) {$c$};  

\draw (\l-\z,0) arc (180:0:{\z/ 2}  and 0.2);
\node[above] at (\l-\z/1.8,0.1) {$z$};   
\node[above] at (\l-\z-\z/20,-0.1) {$c$};

\draw (0,0) arc (-180:0:{\zz/ 2}  and 0.2);
\node[below] at (\zz/1.8,-0.1) {$z'$};   
\node[below] at (\zz+0.1,0) {$c'$}; 

\draw (\u,0) arc (-180:0:{\zz/ 2}  and 0.2);
\node[below] at (\u+\zz/1.8,-0.1) {$z'$};   
\node[below] at (\u+\zz+0.1,0) {$c'$};

\foreach \x in {0,\u,2*\u}
	 \draw[dashed] (\x,0.85) -- (\x ,-0.2);

\end{tikzpicture}
\caption{The case $w=x^{s}, 2\leq s <3, |w|=n=2l+r, |cz|\geq r$.}
    \label{fig:2+b}
\end{figure}
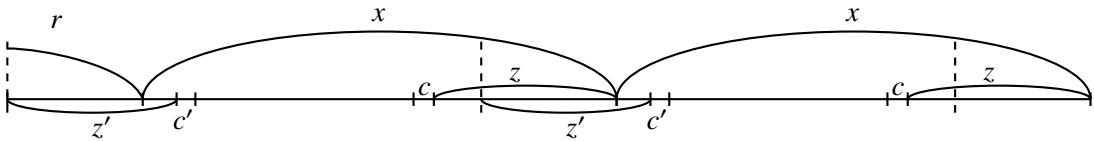

\medskip

Let us count the borders that are suffixes of the rightmost occurrences of closed factors of $w$. 

Any closed factor ending at $j<l+r+1$ is equal to the factor of the same length ending at $j+l$, and thus it is not the rightmost occurrence of this closed factor. The factor $w_{l+r-|z|}\cdots w_{l+r+1} =czw_{l+r+1}$ is unioccurrent in $w$, thus it cannot be a factor of any border. The same is true for $w_{l-1}\cdots w_{l+|z'|+1}=w_{l-1}z'c'$. Thus, if a border starts at the position $i\leq l$, then $i\geq l+r-|z|+1$ (otherwise it contains $czw_{l+r+1}$). This border must end at $j \leq l+|z'|$, otherwise it contains $w_{l-1}z'd'$. 

Summing up, the borders can occur in $w$ as  
$w_{i}\cdots w_{j}$ for $i\in [l+1, n]$ and $j\in [l+r+1, n]$ 
or for $i\in [l+r-|z|+1, l]$ and $j\in[l+r+1, l+|z'|]$ (see Fig.~\ref{fig:2+b}). Without loss of generality, we let $|z|\geq |z'|$. 

There are $\sum_{i=l+1}^{l+r} l + \sum_{i=l+r+1}^n (n-i) $   factors from the first group and $(|z|-r)(|z'|-r)\leq (|z|-r)^2$ factors from the second group. 

Using $r=n-2l$ and $|z|\leq l+r-Cn$, we obtain the following.
\[
\cl'(w) \leq \sum_{i=l+1}^{l+r} l + \sum_{i=l+r+1}^n (n-i) + (|z|-r)^2 \leq lr + \frac{l(l+1)}{2} + (l-Cn)^2 = ln-\frac{l^2}{2}+\frac{l}{2}+Cn(Cn-2l).
\]


Since $\frac{1}{3}\leq C \leq \frac{l}{2}$, for any $l$ we have the maximum of the expression $Cn(Cn-2l)$ when $C=\frac{1}{3}$. 

Thus, we get $\cl'(w) \leq -\frac{l^2}{2} +l\left(\frac{n}{3}+\frac{1}{2}\right) +\frac{n^2}{9}$, and for $C>\frac{1}{3}$, the inequality still holds.
The maximum of this expression is when $l=\frac{n}{3}+\frac{1}{2}$.  Thus, in this case $\cl'(w)\leq \frac{1}{6}n^2 + \frac{n}{6}+\frac{1}{8}.$ Since $\cl'(w)$ is integer, we actually have $\cl'(w)\leq \frac{1}{6}n^2 + \frac{n}{6}.$



3)\ \ Let $w=x^s u= u' x'v'x'$ for non-empty words $u,u',v$. Since $|x|,|x'|> \frac{n}{3}$, the rightmost occurrence of $x'$ must start before position $2l$. In this case, $z$ is the longest repeated prefix of $t$, and $l-|z|=Cn$.

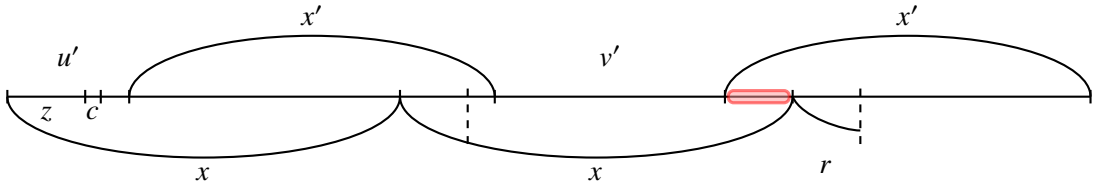
\begin{figure}[ht]
\centering
\begin{tikzpicture}[thick,scale=0.9, every node/.style={scale=0.75}]
\centering
\pgfmathsetmacro\l{16}
\pgfmathsetmacro\u{5.4} 
\pgfmathsetmacro\uu{5.8} 
\pgfmathsetmacro\r{1} 
\pgfmathsetmacro\i{1.8} 
\pgfmathsetmacro\v{\l-\u-\u-\i} 
\pgfmathsetmacro\z{1.15}
\pgfmathsetmacro\zz{1.5}
\Large
\draw (0,0) -- (\l,0);
\draw (\i,0) arc (180:0:{\u/ 2}  and 0.9);
\node[above] at (\i+\u/2,0.9) {$x'$};

\node[above] at (\i+\u+\v/2,0.3) {$v'$};
\node[above] at (\i/2,0.3) {$u'$};

\draw (\l-\u,0) arc (180:0:{\u/ 2}  and 0.9);
\node[above] at (\l-\u/2,0.9) {$x'$};

\draw (0,0) arc (-180:0:{\uu/ 2}  and 0.9);
\node[below] at (\uu/2,-0.9) {$x$};
\draw (\uu,0) arc (-180:0:{\uu/ 2}  and 0.9);
\node[below] at (\uu+\uu/2,-0.9) {$x$};

\draw (2*\uu+\r,-0.5) .. controls (2*\uu+\r-\r/4,-0.5) and (2*\uu+\r/8,-0.3)  .. (2*\uu,0);
\node[below] at (2*\uu+\r/2,-0.8) {$r$};
\draw[dashed] (2*\uu+\r,-0.7) -- (2*\uu+\r,0.2);
\draw[dashed] (\uu+\r,-0.7) -- (\uu+\r,0.2);

\path (\l-\u,0) ++(0.05, 0.6 * 0.15) coordinate (a);
    \path (2*\uu,0) ++(-0.05, -0.6 * 0.15) coordinate (b);
    \draw[rect-color = {red}] (a) rectangle (b);

\foreach \x in {0,\i,\z,\z+\z/5,\i+\u, \l-\u,\l,\uu,2*\uu} 
   \draw(\x,-0.1) -- (\x ,0.1);

\node[below] at (\z/2,0.01) {$z$};   
\node[below] at (\z+\z/10,0.01) {$c$};   

\end{tikzpicture}
\caption{\normalfont{The case $w=u'x'v'x'=x^su$, $2\leq s<3$}}
\label{fig:2su_cases}
\end{figure}

Let us consider the overlapping interval of the occurrence of $x'$ at the position $n-l'+1$ and the occurrence of $x$ at the position $l+1$ in $w$ (see the thick red line on Fig.~\ref{fig:2su_cases}). Its length is $m=2l-n+l'$. If $m\geq|z|$, we have the following. Since $|z|=l-Cn$, we get $2l-n+l'\geq l-Cn$, and thus $l+l'+Cn\geq n$, which is not possible.
If $m<|z|$, then, splitting $w$ into three parts starting at the positions 1, $l+1$ and $n-l'+1$, we get $n=l+(l-m)+l' > l+(l-|z|)+l'=l+Cn+l'>n$, a contradiction.
Thus, this case is not possible.

\bigskip
\noindent

4)\ \ Let $w=u'x'v'x'=xvxu$ for non-empty words $u,v,u',v'$ and $x, x'$ such that $|x|-|z|> \frac{n}{3}$ and $|x'|-|z'|> \frac{n}{3}$. 

As in case 1, we count borders that are suffixes of the rightmost occurrences of closed factors of $w$. For each position~$j$, we consider the longest factor $x''$ ending in $j$ which occurs in $w$ earlier, and consider its longest repeated suffix $z''$. At position $j$, we have $C_j =|x''|-|z''|$ such borders. We let $s$ denote the letter preceding $z''$ in the suffix of $x''$, so that $sz''$ is unioccurrent in $x''$.

Let us suppose there is a factor $x''$ such that $|x''|-|z''|> \frac{n}{3}$. We will now show that in this case both occurrences of $x''$ correspond to factors of the two copies of $x'$ (at the same indices). Let us suppose the converse.

\begin{figure}[ht]
\begin{minipage}{0.5\textwidth}
\centering
\begin{tikzpicture}[thick,scale=0.8, every node/.style={scale=0.8}]
    \def\len{10}
    \def\eps{0.1}
    \draw[thick] (0, 0) -- (\len, 0);
    \draw[thick] (0, -\eps) -- (0, \eps);
    \draw[thick] (\len, -\eps) -- (\len, \eps);
    \def\xval{{0.1, 3.35, 4.2, 7.55, 0.3, 3.8, 6.5, 10, 6.3, 6.1}}
    \def\name{{"$x''$", "$x''$", "$x'$", "$x'$"}}
    \def\up{{1, 1, -1, -1}}

    \foreach \i in {0, 1, ..., 9}{
        \path (\xval[\i], 0) coordinate (n\i);
        \draw[thick] (\xval[\i], -\eps) -- (\xval[\i], \eps);
    }
    \draw[thick] (2.1, -\eps) -- (2.1, \eps);
    \draw[thick] (1.9, -\eps) -- (1.9, \eps);
    \foreach \i in {0, 1}{
        \pgfmathtruncatemacro{\k}{int(2 * \i)}%
        \pgfmathtruncatemacro{\j}{int(2 * \i + 1)}%

        \draw (n\k) arc (180:0:{({\xval[\j] - \xval[\k]}) / 2} and 0.7)
            node[midway, above] {\pgfmathparse{\name[\i]}\pgfmathresult};
    }
    \foreach \i in {2, 3}{
        \pgfmathtruncatemacro{\k}{int(2 * \i)}%
        \pgfmathtruncatemacro{\j}{int(2 * \i + 1)}%

        \draw (n\k) arc (-180:0:{({\xval[\j] - \xval[\k]}) / 2} and 0.7)
            node[midway, below] {\pgfmathparse{\name[\i]}\pgfmathresult};
    }
    \draw (n8) arc (180:0:{({\xval[3] - \xval[8]}) / 2} and 0.15);
    \node[above] at (6.8,0.08) {$z''$}; 
    \path (n9) -- (n8) node[midway, above] {$s$};
    
    \draw (2.1,0) arc (180:0:{({\xval[3] - \xval[8]}) / 2} and 0.15);
    \node[above] at (2.6,0.08) {$z''$}; 
    \node[above] at (2.06,0.02) {$s$};

    \path (n0) ++(0.05, 0.6 * \eps) coordinate (a);
    \path (2.12,0) ++(-0.05, -0.6 * \eps) coordinate (b);
    \draw[rect-color = {red}] (a) rectangle (b);

    \path (n2) ++(0.05,  0.6 * \eps) coordinate (a);
    \path (n8) ++(-0.05, -0.6 * \eps) coordinate (b);
    \draw[rect-color = {red}] (a) rectangle (b);
    
    \path (n6) ++(0.05,  0.6 * \eps) coordinate (a);
    \path (n7) ++(-0.05, -0.6 * \eps) coordinate (b);
    \draw[rect-color = {red}] (a) rectangle (b);
\end{tikzpicture}\subcaption{{\normalfont The rightmost $sz''$ is not a factor of $x'$, no overlap.}}
\end{minipage}
\begin{minipage}{0.5\textwidth}
\centering
\begin{tikzpicture}[thick,scale=0.8, every node/.style={scale=0.8}]
\def\len{10}
    \def\eps{0.1}
    \draw[thick] (0, 0) -- (\len, 0);
    \draw[thick] (0, -\eps) -- (0, \eps);
    \draw[thick] (\len, -\eps) -- (\len, \eps);

    \def\xval{{1, 4.45, 3.9, 7.35, 0.1, 3.6, 6.5, 10, 6.3, 6.1, 3.4, 3.2}}
    \def\name{{"$x''$", "$x''$", "$x'$", "$x'$"}}
    \def\up{{1, 1, -1, -1}}

    \foreach \i in {0, 1, ..., 11}{
        \path (\xval[\i], 0) coordinate (n\i);
        \draw[thick] (\xval[\i], -\eps) -- (\xval[\i], \eps);
    }
    
    \foreach \i in {0, 1}{
        \pgfmathtruncatemacro{\k}{int(2 * \i)}%
        \pgfmathtruncatemacro{\j}{int(2 * \i + 1)}%

        \draw (n\k) arc (180:0:{({\xval[\j] - \xval[\k]}) / 2} and 0.7)
            node[midway, above] {\pgfmathparse{\name[\i]}\pgfmathresult};
    }

    \foreach \i in {2, 3}{
        \pgfmathtruncatemacro{\k}{int(2 * \i)}%
        \pgfmathtruncatemacro{\j}{int(2 * \i + 1)}%

        \draw (n\k) arc (-180:0:{({\xval[\j] - \xval[\k]}) / 2} and 0.7)
            node[midway, below] {\pgfmathparse{\name[\i]}\pgfmathresult};
    }

    \draw (n8) arc (180:0:{({\xval[3] - \xval[8]}) / 2} and 0.15);
    \node[above] at (6.65,0.05) {$z''$}; 

    \path (n9) -- (n8) node[midway, above] {$s$};
    
    \draw (n10) arc (180:0:{({\xval[1] - \xval[10]}) / 2} and 0.15);
    \node[above] at (3.8,0.05) {$z''$}; 

    \path (n11) -- (n10) node[midway, above] {$s$};

    \path (n0) ++(0.05, 0.6 * \eps) coordinate (a);
    \path (n10) ++(-0.05, -0.6 * \eps) coordinate (b);
    \draw[rect-color = {red}] (a) rectangle (b);

    \path (n2) ++(0.05,  0.6 * \eps) coordinate (a);
    \path (n8) ++(-0.05, -0.6 * \eps) coordinate (b);
    \draw[rect-color = {red}] (a) rectangle (b);
    
    \path (n6) ++(0.05,  0.6 * \eps) coordinate (a);
    \path (n7) ++(-0.05, -0.6 * \eps) coordinate (b);
    \draw[rect-color = {red}] (a) rectangle (b);
    
\end{tikzpicture}\subcaption{{\normalfont The rightmost $sz''$ is not a factor of $x'$, there is an overlap.}}
\end{minipage}

\begin{minipage}{0.5\textwidth}
\centering
\begin{tikzpicture}[thick,scale=0.8, every node/.style={scale=0.8}]
    \def\len{10}
    \def\eps{0.1}
    \draw[thick] (0, 0) -- (\len, 0);
    \draw[thick] (0, -\eps) -- (0, \eps);
    \draw[thick] (\len, -\eps) -- (\len, \eps);

    \def\xval{{0.1, 3.35, 6.2, 9.55, 2.55,5.9, 6.5, 10, 2.3, 2.1}}
    \def\name{{"$x''$", "$x''$", "$x'$", "$x'$"}}
    \def\up{{1, 1, -1, -1}}

    \foreach \i in {0, 1, ..., 9}{
        \path (\xval[\i], 0) coordinate (n\i);
        \draw[thick] (\xval[\i], -\eps) -- (\xval[\i], \eps);
    }

    \foreach \i in {0, 1}{
        \pgfmathtruncatemacro{\k}{int(2 * \i)}%
        \pgfmathtruncatemacro{\j}{int(2 * \i + 1)}%

        \draw (n\k) arc (180:0:{({\xval[\j] - \xval[\k]}) / 2} and 0.7)
            node[midway, above] {\pgfmathparse{\name[\i]}\pgfmathresult};
    }

    \foreach \i in {2, 3}{
        \pgfmathtruncatemacro{\k}{int(2 * \i)}%
        \pgfmathtruncatemacro{\j}{int(2 * \i + 1)}%

        \draw (n\k) arc (-180:0:{({\xval[\j] - \xval[\k]}) / 2} and 0.7)
            node[midway, below] {\pgfmathparse{\name[\i]}\pgfmathresult};
    }

    \draw (n8) arc (180:0:{({\xval[1] - \xval[8]}) / 2} and 0.15);
    \node[above] at (2.75,0.05) {$z''$}; 

    \path (n9) -- (n8) node[midway, above] {$s$};

    \path (n0) ++(0.05, 0.6 * \eps) coordinate (a);
    \path (n8) ++(-0.05, -0.6 * \eps) coordinate (b);
    \draw[rect-color = {red}] (a) rectangle (b);

    \path (n4) ++(0.05,  0.6 * \eps) coordinate (a);
    \path (n5) ++(-0.05, -0.6 * \eps) coordinate (b);
    \draw[rect-color = {red}] (a) rectangle (b);
    
    \path (n6) ++(0.05,  0.6 * \eps) coordinate (a);
    \path (n7) ++(-0.05, -0.6 * \eps) coordinate (b);
    \draw[rect-color = {red}] (a) rectangle (b);
\end{tikzpicture}

\subcaption{\centering{\normalfont The rightmost $sz''$ is a factor of $x'$, $x'$ ends before\newline the end of the second occurrence of $sz''$.}}
\end{minipage}
\begin{minipage}{0.5\textwidth}
\centering
\begin{tikzpicture}[thick,scale=0.8, every node/.style={scale=0.8}]
\centering
\def\len{10}
    \def\eps{0.1}
    \def\x{3.5}
    \def\xx{3.35}
    \def\zz{1.05}
    \def\z{1}

    \draw[thick] (0, 0) -- (\len, 0);
    \draw[thick] (0, -\eps) -- (0, \eps);
    \draw[thick] (\len, -\eps) -- (\len, \eps);

    \def\xval{{2.5, 2.5+\xx, 5.1, 5.1+\xx, 0.7, 0.7+\x, 6.5, 6.5+\x, 5.1+\xx-\zz, 5.1+\xx-\zz-0.2, 2.5+\xx-\zz, 2.5+\xx-\zz-0.2, 0.7+\x-\z, 0.7+\x-\z-0.2 , 10-\z,10-\z-0.2}}
    \def\name{{"$x''$", "$x''$", "$x'$", "$x'$"}}
    \def\up{{1, 1, -1, -1}}

    \foreach \i in {0, 1, ..., 15}{
        \path (\xval[\i], 0) coordinate (n\i);
        \draw[thick] (\xval[\i], -\eps) -- (\xval[\i], \eps);
    }

    \foreach \i in {0, 1}{
        \pgfmathtruncatemacro{\k}{int(2 * \i)}%
        \pgfmathtruncatemacro{\j}{int(2 * \i + 1)}%

        \draw (n\k) arc (180:0:{({\xval[\j] - \xval[\k]}) / 2} and 0.7)
            node[midway, above] {\pgfmathparse{\name[\i]}\pgfmathresult};
    }

    \foreach \i in {2, 3}{
        \pgfmathtruncatemacro{\k}{int(2 * \i)}%
        \pgfmathtruncatemacro{\j}{int(2 * \i + 1)}%

        \draw (n\k) arc (-180:0:{({\xval[\j] - \xval[\k]}) / 2} and 0.7)
            node[midway, below] {\pgfmathparse{\name[\i]}\pgfmathresult};
    }

    \draw (n8) arc (-180:0:{({\xval[3] - \xval[8]}) / 2} and 0.1)
            node[midway, below] {$z''$};

    \path (n9) -- (n8) node[midway, below] {$s$};
    
    \draw (n10) arc (-180:0:{({\xval[1] - \xval[10]}) / 2} and 0.1)
            node[midway, below] {$z''$};

    \path (n11) -- (n10) node[midway, below] {$s$};
    
    \draw (n12) arc (180:0:{({\xval[5] - \xval[12]}) / 2} and 0.1);
    \node[above] at (3.7,0.08) {$z'$}; 

    \path (n12) -- (n13) node[midway, above] {$d$};
    
    \draw (n14) arc (180:0:{({10 - \xval[14]}) / 2} and 0.1);
    \node[above] at (10 -0.5 ,0.08) {$z'$}; 

    \path (n14) -- (n15) node[midway, above] {$d$};

  \node[below] at (5.1+\xx ,0) {\textcolor{red}{$\mathbf j$}};  
  \node[below] at (2.5+\xx ,0) {\textcolor{red}{$\mathbf k$}};  
    
\end{tikzpicture}\subcaption{\centering{\normalfont The rightmost $sz''$ is a factor of $x'$, $x'$ ends before\newline the beginning of the second occurrence of $sz''$.}}
\end{minipage}

\caption{\normalfont{Possible locations of ${x''}$ in $w=u'x'v'x'=xvxu$.}}
\label{fig:uxvx_cases}
\end{figure}
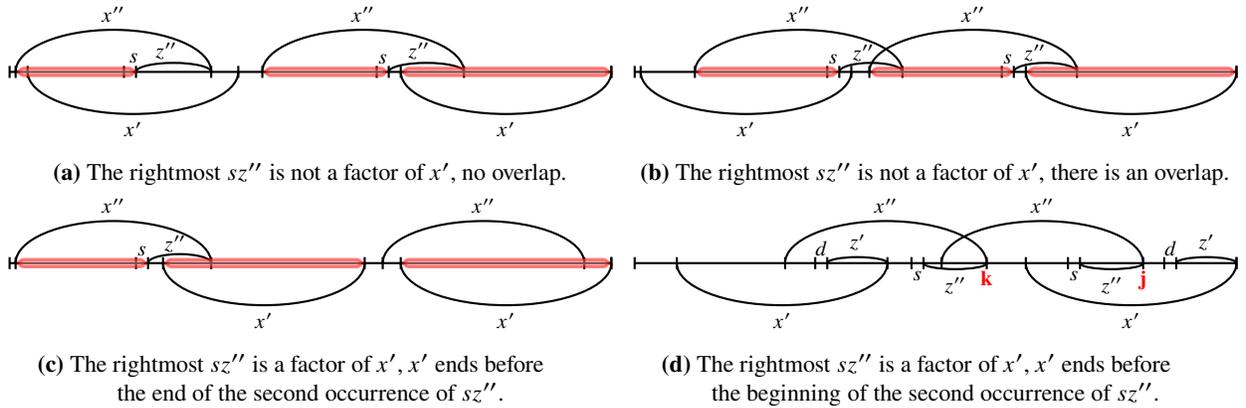

Denoting by $j$ the position of this rightmost occurrence of $x''$, we first consider the case $j-|z''| < n- |x'|$, i.e., the suffix $z''$ of the rightmost occurrence of $x''$ is not a factor of $x'$. 
Then the sum of lengths of the suffix $x'$ of $w$, and the two occurrences of $x''$ without the suffixes $z''$ would be greater than $n$. The mentioned parts of $w$ are drawn with bold on Fig.~\ref{fig:uxvx_cases} (a) (occurrences of $x''$ do not overlap) and (b) (occurrences of $x''$ overlap). Let us note that the overlap cannot be longer than $z''$ due to unioccurrence of $sz''$ in $x''$.

If the suffix $sz''$ of the rightmost occurrence of $x''$ (ending in $j$) is a factor of $x'$, then the previous occurrence of $x''$ (let us say, ending in position $k$) must overlap $x'$ in a way it does not add an extra occurrence of $sz''$ in $x''$. 
If the suffixes $sz''$ of the two occurrences of $x''$ are factors of the two occurrences of $x'$, then they occur at different positions in $x'$ (since we chose $x'$ as the longest repeated prefix). 
So, we have two occurrences of $sz''$ in $x'$ and also in $x''$. 
Thus, the occurrence of $sz''$ ending in $k$ is not a factor of the previous $x'$. 
Hence either $x'$ starts after the beginning of this occurrence of $sz''$, or $x'$ ends before the ending of it. The possibilities are shown on Fig.~\ref{fig:uxvx_cases} (c) and (d). In the first case the sum of lengths of the two occurrences of $x'$ and $|x''|-|z''|$ is greater than $n$. 

In the second case we consider the occurrences of $dz'$.
The first occurrence of $dz'$ as suffix of the first of $x'$ starts at the position $|u'|+|x'|-|z'|$ and ends at the position $n-|x'|-|v'|$.
Its occurrence in the second $x''$ starts at the position $|u'|+|x'|-|z'| + (j-k)$.
In order for it to occur in $x'$ only as suffix the inequality $|u'|+|x'|-|z'| + (j-k) < n- |x'|$, i.e.\ $j-k \leq n-|x'|-(|x'|-|z'|)-|u'|$, must hold. Thus, $j-k < \frac{n}{3}$, hence $C_j < \frac{n}{3}$. Thus, any border that is not a factor of $x$ cannot be a border of more than $\frac{n}{3}$ closed factors. 

\bigskip

Let us count the number of closed factors in $w$. Without loss of generality we assume that $C\leq C'$. We count separately closed factors which are factors of $x'v'x'$, and closed factors starting inside $u$ which do not occur later in the word. 
 By what we proved above, we have at most $Cn$ closed factors starting at each position inside $u$ which do not occur later in the word. For the closed words inside $xvx$, we use the same counting as in the first case. We have the following.
\begin{equation*}
\begin{split}
\cl'(w) &  = |u'|Cn + \sum\limits_{j=|u'|+C'n+1}^{|u'|+2C'n} (j-C'n-|u'|) + \sum\limits_{j=|u'|+2C'n+1}^{n} C'n\\[2ex]
&= |u'|Cn + \frac{C'n(C'n+1)}{2} + C'n(n-|u'|-2C'n) \\[2ex]
 &= -\frac{3}{2} n^2C'^2 + \left(n^2+\frac{n}{2} \right) C' + (C-C')|u'|n\leq -\frac{3}{2}n^2C'^2 + \left(n^2+\frac{n}{2} \right) C'.
\end{split}
\end{equation*}
The last inequality is due to $C\leq C'$. As in case 1), we obtain $\cl'(w) \leq \frac{1}{6}n^2 + \frac{1}{6}n$.

The maximal bound among the obtained ones is $\frac{1}{6}n^2+\frac{1}{6}n$. \quad
$\square$

\bigskip

Now we are ready to prove the main theorem of this section.

We let $\pref_i(w)$ and $\suff_i(w)$ denote the prefix and the suffix of $w$ of length~$i$, respectively. 

\medskip

\textit{Proof }[of Theorem~\ref{th:main}]. 
Let us suppose that $w$ is a word of length $n$ with more than 
$\frac{1}{6}n^2 + \frac{1}{6}n$  long closed factors. 

In other words, $\sum\limits_{j=1}^{n} (\cl'(\pref_{j}(w)) - \cl'(\pref_{j-1}(w)) \geq \frac{1}{6}n^2 + 
\frac{1}{6}n$. It would mean that at least one of the terms in the sum, let us say the $i$-th one, is 
greater than $\frac{i}{3}$. 

Let us consider the largest index $i$ satisfying $\cl'(\pref_{i}(w)) - \cl'(\pref_{i-1}(w)) 
> \frac{i}{3}$, i.e., for all $j>i$, we have $\cl'(\pref_{j}(w)) - \cl'(\pref_{j-1}(w)) < \frac{i}{3}$.
Using Lemma~\ref{upperC}, we can bound the number of distinct long closed factors of $w$ the following way.
\[\cl'(w)  < \frac{1}{6}i^2 + 
\frac{1}{6}i +
\sum\limits_{j=i+1}^{n} \frac{j}{3}  = 
\frac{n^2}{6} +\frac{n}{6}-\frac{i}{6}
\]
Thus, the word $w$ has less than $\frac{n^2}{6}+\frac{n}{6}$ long closed factors. 
$\square$

The following proposition provides a lower bound on $\cl(n)$.
\begin{proposition} \label{pr:lowerbound} 
For each $n$ there exists a word with at least $ \frac{n^2}{6}-\frac{2n}{3}\lceil\log\frac{n}{3}\rceil-\frac{5n}{6}$ distinct closed factors. 
\end{proposition}
\textit{Proof.}
For lengths $n$ divisible by 3 we can e.g.\ take cubes of prefixes of de~Bruijn words. Let us count the borders of closed factors. Any word longer than $\frac{n}{3}$ which occurs in a word at least twice is a border of a closed factor of length longer than  $\frac{2n}{3}$. Due to Proposition~\ref{le:cube}, all such factors are distinct. For each length $j$, $\frac{n}{3}<j\leq \frac{2n}{3}$, there are $\frac{2n}{3}-j+1$ factors of length $j$ which occur twice. Thus, we have the following number of closed factors with borders of lengths longer than  $\frac{n}{3}$: \[\sum\limits_{j=\frac{n}{3}+1}^{\frac{2n}{3}} \left(\frac{2n}{3}-j+1\right)=\frac{n^2}{18}+\frac{n}{6}.\]  

Let us consider factors occurring at the positions $1,2,\dots, \frac{n}{3}- \lceil \log\frac{n}{3}\rceil$ and of length at least $\lceil \log\frac{n}{3}\rceil$ and at most $\frac{n}{3}$. All these $\left(\frac{n}{3}- \lceil \log\frac{n}{3}\rceil \right)\left(\frac{n}{3}- \lceil \log\frac{n}{3}\rceil +1 \right)$ factors are pairwise distinct by construction. Each of them occurs in the word at least twice, and hence is a border of a closed factor.

Thus, there are at least $\frac{n^2}{6} -\frac{2n}{3}\lceil\log\frac{n}{3}\rceil+\frac{n}{6}$ distinct closed factors in this word.

Words of lengths not divisible by 3 can be obtained by shortening a word of next length divisible by 3. By the argument used, for example, in Proposition~\ref{pr:upper}, a suffix of length 1 or 2 can add at most $n$ closed factors. 
$\square$

\medskip

We remark that in the proof of the previous proposition, one can get a slightly better bound by adding small borders (of length smaller than $\lceil \log\frac{n}{3}\rceil$) and counting a bit more carefully lengths non-divisible by 3.

\section{Concluding remarks}



In this paper, we showed an asymptotics of $\sim\frac{n^2}{6}$ for the maximal number of distinct closed factors in a finite word of length $n$. We call such words closed-rich. The upper bound $ \frac{1}{6}n^2+\frac{7}{6}n+1$ has been proved in Theorem \ref{th:main1}, and the lower bound of $\frac{n^2}{6}-\frac{2n}{3}\lceil \log \frac{n}{3}\rceil -\frac{5n}{6}$  has been shown in Proposition \ref{pr:lowerbound}. The question about precise maximal number of distinct closed factors in a closed-rich word remains open even for binary alphabet:

\begin{question}\label{exact}
What is the exact formula for the maximal number of distinct closed factors in a finite rich word?
\end{question}

We remark that linear summand in the upper bound, as well as the constant in $n\log n$ summand in lower bound can be pushed with a more techincal analysis of the same arguments. However, although the bounds can be squeezed to give a smaller interval, the same method does not seem to give the exact value.

Based on numerical experiments, we also conjecture that closed-rich words are cubes or words of exponent close to~3.
Table~\ref{tab:binary} shows the maximal number of closed factors that a binary word of a given length can contain.
\begin{table}[H]
\caption{The maximal number of closed factors for binary words of length $n$.}
    \centering
    \begin{tabular}{c|cccccccccccc}
        $\mathbf n$ & \textbf{1} & \textbf{2} & \textbf{3} & \textbf{4} & \textbf{5} & \textbf{6} & \textbf{7} & \textbf{8} & \textbf{9} & \textbf{10} & \textbf{11} & \textbf{12} \\
    \hline
        $\max_{|w|=n}{\cl(w)}$ & 2 & 3 & 4 & 6 & 8 & 10 & 12 & 15 & 18 & 21 & 25 & 29 \\
    \hline 
    \hline
    $\mathbf n$ & \textbf{13} & \textbf{14} & \textbf{15} & \textbf{16} & \textbf{17} & \textbf{18} & \textbf{19} & \textbf{20} & \textbf{21} & \textbf{22} & \textbf{23} & \textbf{24} \\
    \hline
     $\max_{|w|=n}{\cl(w)}$ & 33 & 37 & 42 & 48 & 54 & 60 & 66 & 72 & 79 & 86 & 93 & 101 
    \end{tabular}
    \label{tab:binary}
\end{table}

\noindent
Similar calculations have been made in~\cite{Fici_low}, but there were some errors. We made corrections for values $n~=~16, 17, \dots,  20$.
For the lengths we computed, closed-rich words are cubes or close to cubes by their structure.
For example, the word $u=(100101)^3$ of length 18 has 60 closed factors (one can easily verify it). The word $u^-$ has 54 closed factors, and the word $u1$ has 66 closed factors. 

We also defined infinite closed-rich words as words for which there exists a constant $C$ such that each factor of length $n$ contains $Cn^2$ distinct closed factors. A question that naturally arises is that of optimizing the constant: 

\begin{question}\label{qu:absolute_rich}
What is the supremum of the constant for infinite rich words?
\end{question}

\bibliography{mainpaper.bbl}

\end{document}